\documentclass[12pt]{article}
\usepackage{amsmath, amsthm, amssymb}
\numberwithin{equation}{section}
\begin{document}

\author{Ajai Choudhry}
\title{Two pairs of biquadrates with equal sums}
\date{}
\maketitle

\begin{abstract}
In this paper we present a new method of solving the classical diophantine equation $A^4+B^4=C^4+D^4$. Two methods of solving this equation, given by Euler, yield parametric solutions given by polynomials of degrees 7 and 13. Several other parametric solutions are now known, and with the exception of one solution of degree 11, all the published solutions are of degrees $6n+1$ for some integer $n$. The method described in this paper yields new parametric solutions of degrees 21, 39 and 75, that is, degrees that are expressible as $6n+3$.
\end{abstract}
\smallskip

\noindent Mathematics Subject Classification 2020: 11D25 
\smallskip

\noindent Keywords: biquadrates; fourth powers.

\section{Introduction}

This paper is concerned with the classical diophantine equation,
\begin{equation}
A^4+B^4=C^4+D^4, \label{biquadeq1}
\end{equation}
that was first solved by Euler in 1772. In fact, Euler gave two ways of solving Eq. \eqref{biquadeq1} leading to parametric solutions given by polynomials of degrees 7 and 13 respectively (see \cite[p.\ 1062]{La}). Dickson \cite[pp.\ 644--647]{Di} mentions various  methods, found subsequently by several mathematicians, of solving the diophantine equation \eqref{biquadeq1}. Lander \cite[pp.\ 1062--1065]{La} applied  geometric methods to obtain parametric solutions of \eqref{biquadeq1}. Zajta \cite{Za} carried out a survey of the important methods of solving Eq. \eqref{biquadeq1} and obtained a new parametric solution given by polynomials of degree 11. Parametric solutions of degrees 13, 19,  25 and 31 have  been published by various authors (\cite{Br}, \cite{Ch}, \cite{La}, \cite{Za}).

Guy \cite[p.\ 212--213]{Gu} states that  ``a method is known for  generating parametric solutions of $a^4+b^4=c^4+d^4$ which will generate all published solutions from the trivial one $(\lambda, 1, \lambda, 1)$; it will only generate solutions of degree $6n+1$." He also mentions that  there exist solutions of even degree but until now no solution of even degree has been found.

Since the existing methods of solving Eq. \eqref{biquadeq1}  generate only solutions of odd degrees, and solutions of even degree are known to exist, there is interest in devising new ways of solving  Eq. \eqref{biquadeq1}  that may yield solutions of even degrees. Accordingly, we tried to find solutions of \eqref{biquadeq1} by a new method which is presented below in Section \ref{new_method}. While the desired solution of even degree could not be found, the method yielded solutions of degrees 21,  39 and 75, all  of these degrees being expressible as $6n+3$. None of the known solutions are of degree  $6n+3$, hence the new solutions obtained in this paper are   interesting although  solutions of even degrees remain elusive. 

\section{A new method of solving Eq. \eqref{biquadeq1}} \label{new_method}

To solve Eq. \eqref{biquadeq1}, we write,
\begin{equation}
\begin{aligned}
A=a_0x^2+a_1x+a_2, \quad B=b_0x^2+b_1x+b_2, \\
C=a_0x^2-a_1x+a_2, \quad D=b_0x^2-b_1x+b_2,
\label{subs1}
\end{aligned}
\end{equation}
where $a_i, b_i, i=0, 1, 2$, and $x$ are arbitrary parameters. With these values,  Eq. \eqref{biquadeq1} reduces, on transposing all terms to the left-hand side and removing the common factor $8x$, to
\begin{multline}
\quad \quad (a_0^3a_1 + b_0^3b_1)x^6 + (3a_0^2a_1a_2 + a_0a_1^3 + 3b_0^2b_1b_2 + b_0b_1^3)x^4\\
 + (3a_0a_1a_2^2 + a_1^3a_2 + 3b_0b_1b_2^2 + b_1^3b_2)x^2 + a_1a_2^3 + b_1b_2^3=0. \quad \quad 
\label{eqred1}
\end{multline}

We will now choose the parameters $a_i, b_i, i=0, 1, 2$, such that the coefficients of $x^6$ and $x^4$ in Eq. \eqref{eqred1} become $0$. Accordingly, we take,
\begin{equation}
\begin{aligned}
a_1&=b_0^3,\quad  &b_1&=-a_0^3,  \\
a_2&=(a_0^8-b_0^8)u/(3a_0b_0^2), \quad &b_2&=(a_0^8-b_0^8)(u-1)/(3a_0^2b_0), 
\end{aligned}
\label{subs2}
\end{equation}
where $u$ is an arbitrary parameter, and now Eq. \eqref{eqred1} reduces to
\begin{equation}
9a_0^2b_0^2((a_0^8 - b_0^8)u + b_0^8)x^2 + (a_0^8 - b_0^8)^2(3u^2 - 3u + 1)=0. \label{eqred2}
\end{equation}
We now write, without any loss of generality, 
\begin{equation}
a_0=b_0t, \quad x=vb_0^2(t^8-1)/(3t((t^8-1)u+1)), \label{subs3}
\end{equation}
when Eq. \eqref{eqred2} reduces to
\begin{equation}
v^2=(3u^2 - 3u + 1)((1-t^8)u-1). \label{eqred3}
\end{equation}
The  birational transformation defined by
\begin{equation}
\begin{aligned}
v&=Y/(3(t^8-1)), \quad &u&=-(X+3)/(3(t^8-1)), \\
X &= 3(-t^8 + 1)u - 3, \quad &Y &= 3v(t^8 - 1),
\end{aligned}
\label{birat1}
\end{equation}
reduces Eq. \eqref{eqred3} to
\begin{equation}
Y^2 = X(X^2 + 3(t^8 + 1)X + 3t^{16} + 3t^8 + 3). \label{ec1}
\end{equation}

Now Eq. \eqref{ec1} may be considered as the Weierstrass model of an elliptic curve over the function field $\mathbb{Q}(t)$. It was found by trial that a rational point $P$ on the elliptic curve \eqref{ec1} is given by
\begin{equation}
P=((t^4 - t^2 + 1)(t^8 - t^4 + 1)/t^2, (t^{18} + t^{12} + t^6 + 1)/t^3).
\end{equation}
It is readily seen that the point $P$ is not a point of finite order, and we can find infinitely many rational points on the elliptic curve \eqref{ec1} using the group law. In fact, the coordinates of the point $2P$ are readily found and are given by
\begin{equation*}
((t^6 - 2t^4 - 2t^2 + 1)^2/(4t^2), (t^{18} - 17t^{12} - 17t^6 + 1)/(8t^3)).
\end{equation*}

Now, using the coordinates of any rational point on the curve \eqref{ec1}  and the relations \eqref{birat1}, \eqref{subs3}, \eqref{subs2} and \eqref{subs1} in that order, we can find a solution of the diophantine equation \eqref{biquadeq1}. While the point $P$ yields a trivial solution of the diophantine Eq. \eqref{biquadeq1}, the point $2P$ gives a nontrivial solution, which on writing $t=a/b$ and clearing denominators, may be expressed in terms of homogeneous polynomials of degree 21 in arbitrary parameters $p$ and $q$  as follows:
\begin{equation}
A=f(p,  q),   B =f(q,  -p),   C=f(p,  -q),   D=f(q,  p),  
\end{equation}
where
\begin{equation}
\begin{aligned}
f(m, n)&=(m - n)(m^2 + mn + n^2)(2m^{18} + 3m^{15}n^3 \\
& \quad \quad + 23m^{12}n^6 + 6m^9n^9 + 8m^6n^{12} - 9m^3n^{15} - n^{18}).
\end{aligned}
\end{equation}
As a numerical example, taking $p=2, q=1$, yields the following solution of Eq. \eqref{biquadeq1}:
\[
A = 5042177, B = 575226, C = 4659327, D = 3638026.
\]

The points $3P$ and $4P$ yield  solutions  of degrees 39 and 75 respectively. Since these solutions are cumbersome to write, we do not give them explicitly. All the computations to find the solutions were performed on the software MAPLE.

\section{Concluding Remarks}

In an effort to find a solution of even degree of the diophantine Eq. \eqref{biquadeq1}, we explored a new method of attacking the problem, and obtained  new parametric solutions of degrees 21, 39 and 75. However, finding a solution of even degree remains an open problem.

\noindent Postal address: Ajai Choudhry, 13/4 A Clay Square, Lucknow - 226001, India

\noindent E-mail: ajaic203@yahoo.com

\end{document}